\documentclass{amsart}%[12pt]
\usepackage{latexsym}
\usepackage{mathabx}
\usepackage{epsfig}
\usepackage{graphicx}
\usepackage{multirow}
\newtheorem{thm}{Theorem}[section]

\newtheorem{lem}[thm]{Lemma}

\newtheorem{prop}[thm]{Proposition}

\makeatletter

\newcommand{\Rmnum}[1]{\expandafter\@slowromancap\romannumeral #1@}
\makeatother
\newenvironment{prf}{{\sl Proof:}}{\medskip \qed}
\newtheorem{defn}[thm]{Definition}

\newcommand{\E}[0]{\mathcal{E}}

\newcommand{\Co}[0]{\mathbb{C}}
\newcommand{\Z}[0]{\mathbb{Z}}
\newcommand{\Q}[0]{\mathbb{Q}}

\newcommand{\fr}[1]{\{#1\}}
\renewcommand{\P}[0]{\mathbb{P}}

\title{the maximal rank of elliptic Delsarte surfaces}
\author{
Bas~Heijne }%{This work was supported by the Netherlands Organization for Scientific Research (NWO). The author would like to thank prof. Jaap Top for his many helpful comments.}
\email{b.l.heijne@rug.nl}
\address{The Johann Bernoulli Institute for Mathematics and Computer Science (JBI), University of Groningen, 
P.O.Box~407, 9700AK Groningen, the Netherlands.}
\thanks{The author would like to thank Jaap Top for several fruitful discussions}
\thanks{This work was supported by
a grant of the Netherlands Organization for Scientific Research (NWO)}
\thanks{Accepted for publication by the Mathematics of Computation of the American Mathematical Society.}


\date{\today}

\begin{document}

\begin{abstract}
Shioda described in his article \cite{shioda1} a method to compute the Lefschetz number of a Delsarte surface.
In one of his examples he uses this method to compute the rank of an elliptic curve over $k(t)$.
In this article we find all elliptic curves over $k(t)$ for which his method is applicable.
For these curves we also compute the maximal Mordell-Weil rank.
\end{abstract}
\maketitle

\section{Introduction}

Shioda described in \cite{shioda1} a method to compute the Lefschetz number of a Delsarte surface.
In one of the examples he used his method to compute the rank over $k(t)$ of the elliptic surface given by:
\[Y^2=X^3+at^nX+bt^m.\]
Here as in the rest of the article $k$ is a algebraically closed field of characteristic 0.
This rank is bounded by $56$ and equal to 56 if and only if $m$ is even and $2m\equiv3n \mod 2^3\cdot 3^2\cdot 5\cdot7$.
Later in \cite{shioda2} Shioda used this method to find an elliptic surface with rank 68 over $k$.
This is the highest rank known for an elliptic surface over $\Co$.

In this article we will briefly describe the method Shioda used.
This method works for all elliptic curves over $k(t)$ that can be defined by a polynomial of the form:
\begin{equation}
f=\sum_{i=0}^3 t^{a_{i0}}X^{a_{i1}}Y^{a_{i2}}.\label{f}
\end{equation}
The first theorem that we will prove is:
\begin{thm}\label{thm1}
Let $f$ be a polynomial over the field $k(t)$. 
Suppose $f$ defines a curve of genus 1 over $k(t)$.
This curve is birational to a curve given by a polynomial $g$ which also has four terms, and moreover its Newton polygon $\Gamma(g)$ is one of the polygons from figure 1.
\end{thm}
After the proof of this result we will compute the maximal rank of the elliptic curves defined by (\ref{f}) for all corresponding Newton polygons.
This will then lead us to the main theorem of this article:
\begin{thm}
Suppose $E/k(t)$ is an elliptic curve defined by (\ref{f}), then $rank(E(k(t)))\leq 68$.
\end{thm}

\begin{figure}
\includegraphics{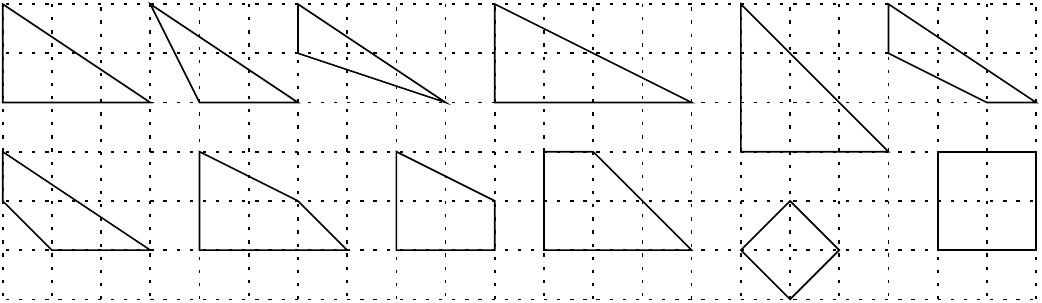}
\caption{All polygons with exactly one interior point and at most 4 corners.}
\end{figure}

\section{Shioda's method/Delsarte Surfaces}
Assume that $f$ is irreducible of the form (\ref{f}), and $f$ defines an elliptic curve $E$ over $k(t)$.

We consider the surface $S\subset \P^3$ over $k$ defined by the homogenized polynomial $F$ corresponding to $f$:
\[F=t^{a_{00}}X^{a_{01}}Y^{a_{02}}Z^{a_{03}}+t^{a_{10}}X^{a_{11}}Y^{a_{12}}Z^{a_{13}}+t^{a_{20}}X^{a_{21}}Y^{a_{22}}Z^{a_{23}}+t^{a_{30}}X^{a_{31}}Y^{a_{32}}Z^{a_{33}}.\]
Let $A_f=(a_{ij})$ be the matrix consisting of the powers appearing in the polynomial.
We assume $A_f$ to be nonsingular.
Then Shioda's method gives us an algorithm to compute the Lefschetz number of the surface.
We will give a slightly adapted description of his method, which appears somewhat more convenient to work with.
For more details we refer to Shioda's orginal paper \cite{shioda1}.

Define $L$ to the subgroup of $(\Q/\Z)^4$ generated by $(1,0,0,-1)A_f^{-1}$,  $(0,1,0,-1)A_f^{-1}$ and  $(0,0,1,-1)A_f^{-1}$.
Define:
\[\Lambda=\left\{(a_i)_i\in L :\begin{array}{l}\forall i: \ a_i\neq0 \mbox{ and }\\ \exists t \in \Z \mbox{ such that }\forall i,\\ ord(ta_i)=ord(a_i) \mbox{ and }\sum_{i=0}^3 \fr{ta_i} \neq 2 \end{array}\right\}.\]
Here $ord$ is the order in the additive group $\Q/\Z$ and $\fr{a_i}$ is the natural bijection between $[0,1)\cap\Q$ and $\Q/\Z$.
\begin{thm}[Shioda]
The Lefschetz number of $S$ is $\lambda=\# \Lambda$.
\end{thm}
\begin{prf}
For the proof see \cite{shioda1}.
\end{prf}

To obtain a formula for the rank of $E/k(t)$ we use the Shioda-Tate formula:
\[rank(E(k(t)))=rank(NS(\E))-\rho_{triv}.\]
Here $\E/k$ is the elliptic surface over $\P^1$ with generic fibre $E/k(t)$.
(So $\E$ and $S$ are birational.)
$NS(\E)$ is the Neron-Severi group of $\E$ and $\rho_{triv}$  is the rank of the subgroup of $NS(\E)$ generated by all components of fibres of $\E\rightarrow \P^1$, together with the $0$-section.
So $\rho_{triv}=2+\sum(m_v-1)$ the sum over all $v\in\P^1$, and $m_v$ is the number of irreducible components of the fibre over $v$.
The rank of $NS(\E)$ is denoted by $\rho$ and can be computed using the formula
\[\rho=h^2-\lambda.\]
Here $h^2$ is the second Betti number of $\E$.
So we have
\[rank(E(k(t)))=h^2-\lambda-\rho_{triv}.\]

\section{Genus calculation.}
%In the previous section we assume that $f$ defines a curve over $k(t)$ of genus 1.
To see which $f$ defines a genus 1 curve, we describe a method which calculates the genus for a given $f$ as in (\ref{f}).
%\[f=t^{a}X^{a_1}Y^{a_2}+t^{b}X^{b_1}Y^{b_2}+t^{c}X^{c_1}Y^{c_2}+t^{d}X^{d_1}Y^{d_2}.\]
To do this we will first need three definitions:
\begin{defn}
An \emph{integral polygon} is the convex hull of a finite subset of $\Z^2$.
\end{defn}

\begin{defn}
Take $f=\sum_{(a,b)\in S}\alpha_{(a,b)} X^{a}Y^{b}$ in the ring of Laurent polynomials $k[X^{\pm 1},Y^{\pm 1}]$, with all  $\alpha_{(a,b)}\not=0$ and $S$ a finite subset of $\Z^2$.
Define the \emph{Newton polygon}, $\Gamma(f)$, of $f$ as the convex hull of $S$.
\end{defn}

\begin{defn}
Take $f=\sum_{(a,b)\in S}\alpha_{(a,b)} X^{a}Y^{b} \in k[X^{\pm 1},Y^{\pm 1}]$.
For every edge, $\gamma$, of the Newton polygon define
$f_{\gamma}=\sum_{(a,b)\in S\cap \gamma} \alpha_{(a,b)} X^{a}Y^{b}$. 
We say that $f$ is nondegenerate with respect to its Newton polygon if for every $\gamma$ we have $f_{\gamma}$, $\frac{\partial f_{\gamma}}{\partial X}$ and $\frac{\partial f_{\gamma}}{\partial Y}$ generate the unit ideal in $k[X^{\pm 1},Y^{\pm 1}]$.
\end{defn}
Note that $f$ is nondegenerate with respect to its Newton polygon if all $f_{\gamma}$ only have simple roots.

We now give the following theorem.
\begin{thm}
Let $f(X,Y)\in k[X^{\pm 1},Y^{\pm 1}]$ be absolutely irreducible.
Then the curve $C$ defined by $f$ has genus
\[g\leq\# \{ \mbox{integral points in the interior of } \Gamma(f) \}.\]
Equality holds if $f$ is nondegenerate with respect to its Newton polygon and the singular points of the projective closure of $C$ in $\P^2$ are all among $(0:0:1)$, $(0:1:0)$ and $(1:0:0)$.
\end{thm}
\begin{prf}
See \cite[Theorem 4.2]{beelen1}
\end{prf}

\begin{lem}
Let $C$ be the projective closure in $\P^2$ of a curve over $k(t)$ defined by a polynomial $f$ as in (\ref{f}). 
Assume that $A_f$ is nonsingular. 
Then $C$ does not have singular points outside the points $(0:0:1)$, $(0:1:0)$ and $(1:0:0)$ over $\bar{k(t)} $.
\end{lem}
\begin{prf}
The curve $C$ is given by the polynomial:
\[\tilde{f}=\sum_{i=0}^3 t^{a_{i0}}X^{a_{i1}}Y^{a_{i2}}Z^{m-a_{i1}-a_{i2}}=0,\]
where $m$ is a positive integer.
Since $A_f$ is nonsingular we can see that:
\[B=
\left(
\begin{array}{cccc}
a_{00}&a_{01}& a_{02}& m-a_{01}-a_{02}\\
a_{10}&a_{11}& a_{12}& m-a_{11}-a_{12}\\
a_{20}&a_{21}& a_{22}& m-a_{21}-a_{22}\\
a_{30}&a_{31}& a_{32}& m-a_{31}-a_{32}\\
\end{array}
\right)
,\]
is likewise nonsingular. 

We claim that there exist $a$, $b_1$, $b_2$, $b_3\in \Z$ with the following property.
With $s\in \bar{k(t)}$ satisfying $s^a=t$, the map $(X:Y:Z)\rightarrow (s^{b_1}X:s^{b_2}Y:s^{b_3}Z)$ defined over $k(s)\supset k(t)$ defines an isomorphism from $C$ to $\tilde{C}$.
Here $\tilde{C}$ is the curve given by:
\[\tilde{f}:= s^nX^{a_{01}}Y^{a_{02}}Z^{m-{a_{01}}-a_{02}}+\sum_{i=1}^3X^{a_{i1}}Y^{a_{i2}}Z^{m-{a_{i1}}-a_{i2}},\]
for some $n\in\Z_{> 0}$.

The proof of this claim is an excercise in linear algebra, as follows:
we look for $\alpha$, $\beta$, $\gamma\in\Q$ such that
\[(t^{\alpha}X)^{a_{i1}}(t^{\beta}Y)^{a_{i2}}(t^{\gamma}Z)^{m-a_{i1}-a_{i2}}=t^{a_{i0}}X^{a_{i1}}Y^{a_{i2}}Z^{m-a_{i1}-a_{i2}},\]
for $i=1,2,3$.
This means that $v=(\alpha,\beta,\gamma)^T$ is a solution of 
\[B_{1,1}v=\left(
\begin{array}{c}
a_{10}\\
a_{20}\\
a_{30}\\
\end{array}
\right),\]
in which $B_{1,1}$ is the submatrix of $B$ obtained by deleting the first row and the first column.
Since by assumption $\det (B_{1,1})\neq 0$, the solution $v$ exists and is unique.

Let $B_1$ the the first column of $B$ and $B'$ the matrix obtained by deleting the first column of $B$.
Then since $B$ is nonsingular we find $B'v\neq B_1$.
This means precisely that $n\neq 0$.
Now by taking $a$ a multiple of the denominators of $v$ and of correct sign we can ensure that $n$ is a positive integer.

To prove the lemma, it suffices to prove that $\tilde{C}$ has no singular points outside $(1:0:0)$, $(0:1:0)$ and $(0:0:1) $.
Assume $(x:y:z)$ is a singular point on $\tilde{C}$.
Write $w_0=s^nx^{a_{01}}y^{a_{02}}z^{m-{a_{01}}-a_{02}}$ and $w_i=x^{a_{i1}}y^{a_{i2}}z^{m-{a_{i1}}-a_{i2}}$ for $i \in \{1,2,3\}
$ and $w=(w_0,w_1,w_2,w_3)$.

Note that $wB'=0$. 
Now there are two possibilities. 
The first possibility is that $w_0=0$. Since $\det(B_{1,1})\neq0$ this implies $w=(0,0,0,0)$. This implies that $(x:y:z)$ is one of $(1:0:0)$, $(0:1:0)$ or $(0:0:1)$.

If $w_0\neq0$ we can assume $w_0=1$. Since $\det(B_{1,1})\neq 0$ we now find that $w\in \Q^4$.
Again using that the last three rows of $B$ are linearly independent this gives that $x,y,z$ are algebraic over $\Q$. 
This contradicts the fact that $v_0=1$ and $n>0$.

\end{prf}

\begin{lem}
Let $f$ be as in (\ref{f}).
Assume that $\det(A_f)\neq 0$, then $f$ is nondegenerate with respect to its Newton polygon.
\end{lem}

\begin{prf}
For any edge $\gamma$ we find $f_{\gamma}$ has either two or three terms. 
Four terms on one edge is not possible, since then $\det(A_f)=0$.
The case where $f_{\gamma}$ has only two terms is simple, since $\mbox{char}(k)=0$.
%Then $f_\gamma=\alpha_1X^{a_1}Y^{a_2}+\alpha_2X^{b_1}Y^{b_2}$.
%It is now trivial to find that there are no $X,Y\neq 0$ with  $f_\gamma=X\frac{\partial f_\gamma}{\partial X}=Y\frac{\partial f_\gamma}{\partial Y}=0$, unless $a_1=b_1$ and $a_2=b_2$.

We will only do the case where $f_\gamma$ has three terms.
Without loss of generality we can assume that 
$f_\gamma=X^{a}+Y^{b}+s^nX^{\lambda a}Y^{(1-\lambda)b}$, where $s$ a root of $t$.
Here $n$ is nonzero, since otherwise $\det(A_f)=0$.
Define $\eta=X^{a}/Y^b$. 
Now we assume that $a\neq 0$, then $\frac{\partial f_\gamma}{\partial X}=0$ and $f_{\gamma}-\frac{X}{a}\frac{\partial f_\gamma}{\partial X}=0$ gives:
\[\eta+\lambda s^n \eta^{\lambda}=0\]
\[1+(1-\lambda) s^n \eta^{\lambda}=0\]
%\[\xi^{\lambda(1-\lambda)}=(-\lambda s)^{\lambda} \]
%\[\xi^{\lambda(1-\lambda)}=((\lambda-1) s)^{\lambda-1} \]
Since $s$ is transcendental over $k$ this has no solution.
\end{prf}

Combining the previous two lemma's with the result on the genus gives the following theorem.
\begin{thm}
Let $C$ be a curve over the field $k(t)$, defined by an absolutely irreducible polynomial of the form given in (\ref{f}).
%Also assume that $f$ is not isomorphic to a curve defined over $k$
Also assume that $\det(A_f)\neq0$ , then the genus of $C$ equals the number of interior lattice points of its Newton polygon.
\end{thm}

\section{The forms of the Newton polygon}
In this section we will define the concept of equivalent polygons. 
We will use this to give a result on equivalence of elliptic curves.

\begin{defn}
We call two polygons $A,B$ \emph{integrally equivalent} if $B$ is the image of $A$ under a linear map given by a matrix in $GL_2(\Z)$ possibly composed with a translation.
%there exists a matrix $M\in GL_2(\Z)$, such that
%$MA$ is a shift of $B$. Here we view elements of $A$ as column vectors.
%(i.e. there exists a linear map with integer coefficients mapping $A$ to $B$ and vice versa.) 
\end{defn}
Note that being integrally equivalent is an equivalence relation.

If $f\in k[X,Y]$ is irreducible and $f\not\in k\cdot X\cup k\cdot Y$, then $\Gamma(f)$ has at least one point on both the $x$ and $y$-axis.
Furthermore $\Gamma(f)$ is contained in the first quadrant.
Any integral polygon can be shifted in a unique way such that it satisfies these criteria, i.e. it is contained in the first quadrant and it has a point on each of the axis.
We shall consider this to be the \emph{default position} of the polygon. 

\begin{prop}
Let $f(X,Y)=\sum_{(a,b)\in S}\alpha_{(a,b)} X^{a}Y^{b}$ be a bivariate polynomial over a field, defining an irreducible curve $C$. 
Assume that all $\alpha_{(a,b)}\neq 0$.
Given a polygon $A$, in default position,  integrally equivalent to the polygon $\Gamma(f)$, then there exists an irreducible bivariate polynomial $g(X,Y)$, such that $A=\Gamma(g)$.
Moreover the coefficients of $f$ and $g$ will be the same and the curves defined by $f$ and $g$ will be birationally equivalent.
\end{prop}

\begin{prf}
Let $M=\left(\begin{array}{cc}k&l\\m&n\end{array}\right)$ be the matrix such that $M\Gamma(f)$ is a shift of $A$.
Define $g(U,V)=U^{\lambda}V^{\mu}\sum_{(a,b)\in S}\alpha_{(a,b)} U^{ak+bl}V^{am+bn}$.
Here $\lambda$ and $\mu$ are so that $\Gamma(g)$ is in default position.
By definition $g$ has the same nonzero coefficients as $f$.
The birational equivalence between the curves given by $g$ and $f$ is defined by:
\[\phi: Z(g) \longrightarrow Z(f)\]
\[(U,V) \longrightarrow (U^kV^m,U^lV^n).\]
\end{prf}

It is a well known result (see \cite{beelen1}  and \cite{rabinowitz}) that up to integral equivalence there are exactly 16 polygons with exactly one interior point. %ref
Four of these polygons have more than 4 corners. This gives the final result \ref{thm1}.

%\begin{thm}
%Let $f$ be a polynomial with four terms over the field $k(t)$. Let $f$ define a curve of genus 1.
%Then $f$ is birational to a curve given by a polynomial $g$ wich also has four term.
%Moreover $\Gamma(g)$ is one of the polyfons from figure 1.
%\end{thm}
%\includegraphics{met3}

\section{An example}

In this section we present one example of a computation of the rank.
For the other families of elliptic curve we will provide less details. 
The methods employed will remain the same however.

We will consider the elliptic curves over $k(t)$ that are defined by a polynomial of the form:
\[f=t^a+(t^b+t^c)X^3+t^dY^2=0,\]
where $a,b,c,d$ are integers $\geq$ 0 with $c>b$.
We want to find the maximal rank that occurs in this family.

Let $E$ be the curve defined by $f$ and $E'$ the curve defined by:
\[t^{6a}+(t^{6b}+t^{6c})X^3+t^{6d}Y^2=0.\]
Then we have a natural monomorphism $\phi: E(k(t))\longrightarrow E'(k(t))$, defined by $\phi(x(t),y(t))=(x(t^6),y(t^6))$. In particular we find the rank of $E(k(t))$ is at most the rank of $E'(k(t))$. 
So we will restrict ourselves to computing the rank of $E'$.

Divide the equation of $E'$ by $t^{6a}$ and define $\xi=t^{2(b-a)}X$, $\eta=t^{3(d-a)}Y$ and $n=6(c-b)$. Then we see that $E'$ is isomorphic to the curve $\tilde{E}$ defined by:
\[\tilde{f}=1+(1+t^n)\xi^3+\eta^2=0.\]
So we only have to determine the maximal rank for curves of this form.
Note that for all $m>0$ there is an injective map:
\[\tilde{E}(k(t))\longrightarrow \tilde{E}'(k(t)),\]
\[(\xi(t),\eta(t))\longrightarrow (\xi(t^m),\eta(t^m)).\]
Here $\tilde{E}'$ is the curve given by:
\[1+(1+t^{nm})\xi^3+\eta^2=0.\]
From this we see that without loss of generality we can assume that $m|n$.

We will compute the Lefschetz number using the technique from Shioda.
To do this we first homogenize $\tilde{f}$.
This gives:
\[\tilde{F}=Z^{n+3}+T^nX^3+X^3Z^{n}+Y^2Z^{n+1}.\]
Then we compute the matrices $A$ and $A^{-1}$.
\[A=\left(
\begin{array}{cccc}
0&0&n+3&0\\
3&0&0&n\\
3&0&n&0\\
0&2&n+1&0\\
\end{array}
\right) \;\;\; \mbox{ and }\;\;\; 
A^{-1}=\left(\begin{array}{cccc}
-\frac{n}{3(n+3)}& 0 & \frac{1}{3} & 0\\
-\frac{n+1}{2(n+3)}& 0 & 0 &\frac{1}{2}\\
\frac{1}{n+3}& 0 & 0&0\\
\frac{1}{n+3} &\frac{1}{n} &-\frac{1}{n} &0\\
\end{array}
\right)
.
\]
By definition $L$ is the supgroup of $(\Q/\Z)^*$ generated by  
\[w_1=(1,0,0,-1)A^{-1}=(-\frac{1}{3},-\frac{1}{n},\frac{n+3}{3n},0),\]
\[w_2=(0,1,0,-1)A^{-1}=(-\frac{1}{2},-\frac{1}{n},\frac{1}{n},\frac{1}{2}),\]
\[w_3=(0,0,1,-1)A^{-1}=(0,-\frac{1}{n},\frac{1}{n},0).\]
By inspecting these generators we see that $L$ is also generated by:
\[v_1=w_1-w_3=(-\frac{1}{3},0,\frac{1}{3},0),\]
\[v_2=w_2-w_3=(-\frac{1}{2},0,0,\frac{1}{2}),\]
\[v_3=w_3=(0,-\frac{1}{n},\frac{1}{n},0).\]

We see that $L$ consists of elements of the form $iv_3$, $v_1+iv_3$, $2v_1+iv_3$, 
$v_2+iv_3$, $v_1+v_2+iv_3$ and $2v_1+v_2+iv_3$. For each form there are exactly $n$ elements. 
To compute $\lambda$ we have to find out which of these elements lie in $\Lambda$.

Elements of the form $iv_3,v_1+iv_3$ and $2v_1+iv_3$ do not lie in $\Lambda$, since they all have zero as their last coordinate.

An element of the form $v_2+iv_3$ does not lie in $\Lambda $. If $i=0$ this follows from the fact that the second and third coordinate are zero. If $i\neq 0$ then this follows from the the fact that we can compute for all $t$ with $(t,2n)=1$:
\[\fr{\frac{ti}{n}}+\fr{-\frac{ti}{n}}+\fr{\frac{t}{2}}+\fr{-\frac{t}{2}}=2.\]

We will now determine when $v_1+v_2+iv_3\in \Lambda$.
Take $j,m\in \Z_{\geq 0}$  such that $j/m=i/n$ and $(j,m)=1$ and write $v_1+v_2+iv_3=
(\frac{1}{6},-\frac{j}{m},\frac{1}{3}+\frac{j}{m},\frac{1}{2})$.
The condition $\fr{\frac{t}{6}}\neq 0$, $\fr{-\frac{jt}{m}}\neq0$, $\fr{\frac{t}{3}+\frac{jt}{m}}\neq 0$ and $\fr{\frac{t}{2}}\neq0$ is satisfied precisely when $j\neq 0$ and $\frac{j}{m}\neq \frac{2}{3}$. 

In all other cases we have $v_1+v_2+iv_3\in \Lambda$ if and only if there exists a $t$ such that $(t,6m)=1$ and
\[\fr{\frac{t}{6}}+\fr{-\frac{jt}{m}}+\fr{\frac{t}{3}+\frac{jt}{m}}+\fr{\frac{t}{2}}\neq2.\]
It is easy to compute, if $j\neq 0$ and $\frac{j}{m}\neq\frac{2}{3}$:
\[\fr{\frac{t}{6}}+\fr{-\frac{jt}{m}}+\fr{\frac{t}{3}+\frac{jt}{m}}+\fr{\frac{t}{2}}=\left\{
\begin{array}{ll}
1& \mbox{ if } t\equiv 1 \mod 6 \mbox{ and }\fr{\frac{tj}{m}}>\frac{2}{3}\\
2& \mbox{ if } t\equiv 1 \mod 6 \mbox{ and }\fr{\frac{tj}{m}}<\frac{2}{3}\\
3& \mbox{ if } t\equiv 5 \mod 6 \mbox{ and }\fr{\frac{tj}{m}}<\frac{1}{3}\\
2& \mbox{ if } t\equiv 5 \mod 6 \mbox{ and }\fr{\frac{tj}{m}}>\frac{1}{3}\\
\end{array}
\right.
\]
By considering a pair $\pm t$, this means that $v_1+v_2+iv_3\in \Lambda$ if and only if $\fr{\frac{tj}{m}}<\frac{1}{3}$ for some $t\equiv 5\mod 6$, with $(t,6m)=1$.
We now distinguish between the various possibilities :
\begin{itemize}

\item The case $m\leq 3$ is easy and leads to $(v_1+v_2+iv_3)\not\in \Lambda$.
This happens precisely when $i\in\{0,n/2,n/3,2n/3\}$.

\item Assume $m>3$ and $3\notdivides m$ or $j\equiv 2\mod 3$.  Then $t\in \Z$ exists with $t\equiv 5\mod 6$ and $t\equiv j^{-1} \mod m$. For this $t$ we find $\fr{\frac{tj}{m}}<\frac{1}{3}$, hence 
$(v_1+v_2+iv_3)\in \Lambda$.

\item In the case that $m>3$, $3|m$, $j\equiv 1\mod 3$, assume moreover that there exists a $c\equiv 2\mod 3$, with $(c,m)=1$ and $\fr{\frac{c}{m}}<\frac{1}{3}$. We can find $t\equiv 5 \mod 6$ such that $t\equiv cj^{-1} \mod m$. For that $t$ we have $\fr{\frac{tj}{m}}<\frac{1}{3}$. This means $(v_1+v_2+iv_3)\in \Lambda$. This happens for all $m>3$ except when $m \in \{6,12, 30 \}$, as is shown in lemma \ref{lem1} below.  

\item The final case is $m>3$, $3|m$, $j\equiv 1\mod 3$ and there exists no $c\equiv 2\mod 3$, with $(c,m)=1$ and $\fr{\frac{c}{m}}<\frac{1}{3}$. Assume that $v_1+v_2+iv_3\in \Lambda$. 
Then $t\equiv 5 \mod 6$ exists, coprime to $6m$ such that $\fr{\frac{tj}{m}}<\frac{1}{3}$. Hence $c=jt$ satisfies $c\equiv 2\mod 3$, $gcd(c,m)=1$ and $\fr{\frac{c}{m}}<\frac{1}{3}$, contrary to our assumption.

So in this case we find $(v_1+v_2+iv_3)\not\in \Lambda$. By the following lemma, this final possibility for $m$ and $j$ happens only if $m\in \{6,12, 30 \}$.
In other words only if $i\in\{\frac{n}{6}, \frac{n}{12}, \frac{7n}{12}, \frac{n}{30}, \frac{7n}{30},\frac{13n}{30},\frac{19n}{30}\}$. 
\end{itemize}

\begin{lem}\label{lem1}% We only use a weaker result here in another part we use a stronger result.
6, 12 and 30 are the only integers $n>3$ with the property that there does not exist a prime $p\equiv 2\mod 3$ such that $3p<n$ and $p\not|n$.% satisfy $n\in\{6,12,30\}$.
\end{lem}
\begin{prf}
If $n$ satisfies this property then it can be written as $n=Kp_1p_2\ldots p_t$, with the $p_i$ all primes with $p_i\equiv 2 \mod 3$ and $3p_i<n$.
Order the $p_i$ such that $p_i<p_{i+1}$.
We construct the number $N=3p_1\ldots p_{t-1}+p_t$ and see that it has a prime $p\equiv 2\mod 3 $ dividing it, with $p\not=p_i$.
If $n> 51$ we find:
\[p/n\leq N/n = \frac{3}{Kp_t}+\frac{1}{Kp_1\ldots p_{t-1}}\leq\frac{3}{17}+\frac{1}{2\cdot 5\cdot 11}<\frac{1}{3}.\]
This means $3p<n$, but $p$ is not any of the $p_i$, a contradiction. 
So if $n$ satisfies the conditions of the lemma we have $n\leq 51$. Checking the lemma for $n\leq 51$ is easy.
\end{prf}
\\
The cases $v_1+v_2+iv_3$ and $2v_1+v_2+iv_3$ are similar, since $-(v_1+v_2+iv_3)=2v_1+v_2+(n-i)v_3$ and the fact that $v\in \Lambda \Leftrightarrow -v\in \Lambda$. %Dit kan later mooier gemaakt worden.

To ensure that all the special values $i\in\{0,\frac{n}{2},\frac{n}{3},\frac{2n}{3}, \frac{n}{6}, \frac{n}{12}, \frac{7n}{12}, \frac{n}{30}, \frac{7n}{30},\frac{13n}{30},\frac{19n}{30}\}$ encountered in the calculations are actually integers we assume that $60|n$. %We can do that without losing our maximal solution, since if the rank is maximal for $n$ it is also maximal for $60n$.
In that case we find $\lambda=2n-22$.

To compute the rank of the curve we still have to compute both the $\rho_{triv}$ and $ h^2$.
Both of these we will compute for the curve in short Weierstrass form.
Define $\tilde{\eta}=(1+t^n)\eta$ and $\tilde{\xi}=(1+t^n)\xi$ then we get the formula:
\[\tilde{\eta}^2+\tilde{\xi}^3+(1+t^n)^2=0.\]
Here we use theory explained in \cite{schuett} to show that the second Betti number is $h^{2}=4n-2$.
We also compute:
\[\Delta=-432(t^{n}+1)^4.\]
\[j=0.\]
From this we see that the elliptic surface has $n$ singular fibres of type \Rmnum{4} at the roots of $t^{n}+1=0$ and no other singular fibres. So we find $\rho_{triv}=(2n+2)$. 
For details on the roman numeral notation see \cite{sil1} or \cite{schuett}.

Combining these facts gives:
\[r=h^{2}-\lambda-\rho_{triv}=4n-2-(2n-22)-(2n+2)=18.\]
This concludes the example and we find that the rank of $E$ over $k(t)$ is $\leq 18$ and it equals $18$ when $60|n$.

\section{Results}
The following table is a complete list of all integer polygons  with exactly one interior point and at most four corners, up to equivalence. 
For each polygon we give a list of curves over $k(t)$, such that any elliptic curve with Newton polygon equal to the given one can be injected in one of these curves.
As a consequence we create 42 families of elliptic curves over $k(t)$. Any other Delsarte elliptic curve over $k(t)$ can be injected into at least one of the 42 families. 
%Dit moet mooier.

Note that there is a choice in the contruction of this table.
In most cases the term $t^n$ could be placed in front of another term of the formula defining the curve. 
The choices made are partially arbitrary and have partially been motivated by ease of computing the rank.

The computations done to fill this table can be found in the two appendices.
\\

\begin{tabular}{|c|c|c|c|c|}
\hline
Picture& Name &Form with maximal rank& Maximal rank&Occurring for $n$\\
\hline
\multirow{7}{*}{\includegraphics{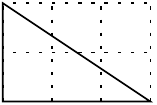}}
&$E^{1a}_n$&$1+t^n+X^3+Y^2$ & 68&360\\%
&$E^{1b}_n$&$1+t^nX+X^3+Y^2$&56&840\\%
&$E^{1c}_n$&$1+t^nX^2+X^3+Y^2$&9&20\\%
&$E^{1d}_n$&$1+(1+t^n)X^3+Y^2$&18&60\\%
&$E^{1e}_n$&$t^n+Y+X^3+Y^2$&68&360\\        %Note:1a and 1e isomorphic.
&$E^{1f}_n$&$1+t^nXY+X^3+Y^2$&9&10\\        %Note:1c and 1f isomporphic.
&$E^{1g}_n$&$1+X^3+(1+t^n)Y^2$&4&6 \\%
\hline

\end{tabular}

\begin{tabular}{|c|c|c|c|c|}
\hline
Picture& Name &Form with maximal rank& Maximal rank&Occurring for $n$\\
\hline
\multirow{5}{*}{\includegraphics{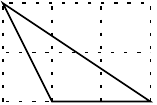}} 
&$E^{2a}_n$& $(1+t^n)X+X^3+Y^2$ & 24& 24\\
&$E^{2b}_n$&$t^nX+X^2+X^3+Y^2$&3&12 \\%
&$E^{2c}_n$&$X+(1+t^n)X^3+Y^2$&24&24 \\ %birational to 2a
&$E^{2d}_n$&$t^nX+XY+X^3+Y^2$&3&12 \\%isomorphic to 2b
&$E^{2e}_n$&$X+X^3+(1+t^n)Y^2$&6&12 \\%
\hline

\multirow{4}{*}{\includegraphics{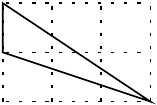}}
&$E^{3a}_n$&$(1+t^n)Y+X^3+Y^2$&18&60\\ %birational to 1d 3b and 3c
&$E^{3b}_n$&$Y+(1+t^n)X^3+Y^2$&18&60\\%birational to 1d 3b and 3c
&$E^{3c}_n$&$Y+X^3+(1+t^n)Y^2$&18&60\\%birational to 1d 3b and 3c
&$E^{3d}_n$&$Y+t^nXY+X^3+Y^2 $ & 1&2\\%
\hline
\multirow{9}{*}{\includegraphics{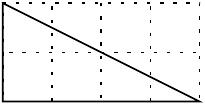}}
&$E^{4a}_n$&$1+t^n+X^4+Y^2$&24&    24 \\ %Birational to 2a and 2c
&$E^{4b}_n$&$1+t^nX+X^4+Y^2$&56&  840 \\%Birational to 1b and t^{1120}+X+X^4+Y^2
&$E^{4c}_n$&$t^n+X^2+X^4+Y^2$&3&   12 \\%check
&$E^{4d}_n$&$1+t^nX^3+X^4+Y^2$&56&840 \\%Birational to 1b\\
&$E^{4e}_n$&$1+(t^n+1)X^4+Y^2$&24& 24 \\ %Birational to 2a and 2c\\
&$E^{4f}_n$&$1+t^nY+X^4+Y^2$&24&   12 \\ %Birational to 2a and 2c\\
&$E^{4g}_n$&$t^n+XY+X^4+Y^2$&3&    12 \\%birational to 4c
&$E^{4h}_n$&$1+t^nX^2Y+X^4+Y^2$&24&12 \\ %Birational to 2a and 2c\\
&$E^{4i}_n$&$1+X^4+(1+t^n)Y^2$&6&  12 \\ %birational to 2e. 
\hline
\multirow{10}{*}{\includegraphics{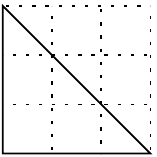}}
&$E^{5a}_n$&$1+t^n+X^3+Y^3$&18&     60\\ %Birational to 1d.
&$E^{5b}_n$&$1+t^nX+X^3+Y^3$&68&   120\\ %Birational to 1a.
&$E^{5c}_n$&$1+t^nX^2+X^3+Y^3$&68& 120\\ %Birational to 1a.
&$E^{5d}_n$&$1+(1+t^n)X^3+Y^3$&18&  60\\ %Birational to 1d.
&$E^{5e}_n$&$1+t^nY+X^3+Y^3$&68&   120\\ %Birational to 1a.
&$E^{5f}_n$&$1+t^nXY+X^3+Y^3$&1&     2\\%
&$E^{5g}_n$&$1+t^nX^2Y+X^3+Y^3$&68&120\\ %Birational to 1a
&$E^{5h}_n$&$1+t^nY^2+X^3+Y^3$&68& 120\\ %Birational to 1a
&$E^{5i}_n$&$1+t^nXY^2+X^3+Y^3$&68&120\\ %Birational to 1a
&$E^{5j}_n$&$1+X^3+(1+t^n)Y^3$&18&  60\\ %Birational to 1d.\\
\hline
\multirow{3}{*}{\includegraphics{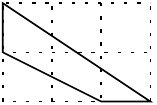}}
&&&&\\
&$E^{6}_n$&$t^nX^2+Y+X^3+Y^2$&9&20\\ %This used to be 9. 
&&&&\\
\hline
\multirow{3}{*}{\includegraphics{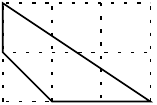}}&&&&\\
&$E^{7}_n$&$t^nX+Y+X^3+Y^2$&56&840\\%birational to 1b
&&&&\\
\hline
\multirow{3}{*}{\includegraphics{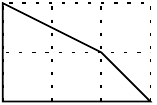}}&&&&\\
&$E^{8}_n$&$1+t^nX^2Y+X^3+Y^2$&56&560\\ %birational to 1b Check the number 420 later used to be 560
&&&&\\
\hline
\multirow{3}{*}{\includegraphics{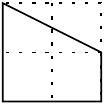}}&&&&\\
&$E^{9}_n$&$t^n+X^2Y+X^2+Y^2$ &3&12\\ %birational to 4c
&&&&\\
\hline
\multirow{3}{*}{\includegraphics{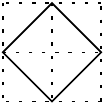}}&&&&\\
&$E^{10}_n$&$t^nX+Y+X^2Y+XY^2$&0&1\\
&&&&\\
\hline
\multirow{3}{*}{\includegraphics{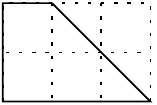}}&&&&\\
&$E^{11}_n$&$t^n+XY^2+X^3+Y^2$&18&120\\
&&&&\\
\hline
\multirow{3}{*}{\includegraphics{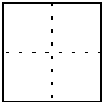}}&&&&\\
&$E^{12}_n$&$t^n+X^2+Y^2+X^2Y^2$&0&1\\
&&&&\\
\hline

\end{tabular}
\appendix
\section{equivalence}

Before we proceed to compute the rank of the elliptic curves in the above table, we give some criteria explaining why the ranks of certain elliptic curves are the same.

First of all if there is an isogeny between two curves $E_1$ and $E_2$ then the ranks of the two curves will be the same.

Secondly if two elliptic curves $E_1$ and $E_2$ are defined over a field $k(t)$ and there is a $k$-linear automorphism $\phi:k(t)\rightarrow k(t)$ bijecting the points of $E_1$ to $E_2$ then these two curves will have the same rank.

Of course sometimes a combination of these two methods can be used to show that two curves have the same rank.
In this appendix we will use these methods to show why certain families of elliptic curves in our table have the same maximal rank.

\begin{defn}
As a matter of terminology, we will say that two curves are $k$-\emph{isogenous} (respectively $k$-\emph{equivalent}) if they are isogenous (respectively isomorphic) after a $k$-linear automorphism.
\end{defn}

%\begin{defn}
%Given two curves $E_1$ and $E_2$ over a rational function field $k(t)$, then the curves are \emph{$k$-equivalent} if there exists a curve $C'$ isomorphic to $C_1$ and a $k$ invariant field automorphism $\phi:k(t)\rightarrow k(t)$ such that $\phi$ bijects $C'$ in $C_2$. Here $\phi$ acts on as a map from $C'$ to $C_2$ by $\phi:(x,y)\rightarrow (\phi(x),\phi(y))$.
%\end{defn}

%\begin{remark}
%This curves $C_1$ and $C_2$ over $k(t)$ are $k$-equivalent if and only if the surfaces over $k$ defined by the same equations are isomorphic. 
%From this it also follows that $k$-equivalence is an equivalence relation. 
%\end{remark}

%In this section we will find some of the $k$-equivalences that occur in the families of curves in our table.

\subsection{1a}
Here we explain the relation between the families of curves $E^{1a}_n$, $E^{1e}_n$, $E^{5b}_n$, $E^{5c}_n$, $E^{5e}_n$, $E^{5g}_n$, $E^{5h}_n$ and $E^{5i}_n$.

Permuting homogeneous coordinates $X$, $Y$, $Z$ gives isomorphisms between the families described by $E^{5b}_n$, $E^{5c}_n$, $E^{5e}_n$, $E^{5g}_n$, $E^{5h}_n$ and $E^{5i}_n$.

A short Weierstrass form for the curves $E^{1e}_n$ is %using the coordinates $\eta=\sqrt{-4}(Y+\frac{1}{2})$ and $\xi=\sqrt[3]{-4}X$.
\[1-4t^n+\xi^3+\eta^2=0.\]
The field automorphism defined by $t\rightarrow \sqrt[n]{-4}t$ brings this precisely in the form described by $E^{1a}_n$.
We conclude that the curves $E^{1a}_n$ and $E^{1e}_n$ are $k$-equivalent.

Using the isomorphism given by $V=\sqrt{3}\frac{X^2-Y^2}{4}+\frac{t^n}{4\sqrt{3}}$ and $U=\frac{X+Y}{2}$ we see that the curve $E^{5b}_n$ is isomorphic to the curve given by:
\[V^2+U^4+\frac{1}{2}t^nU^2+\frac{1}{2}U-\frac{1}{48}t^{2n}=0.\]
Using ideas described in \cite[Ch. 8 ]{cassels} one finds that this curve is isomorphic to the curve given by:
\[\eta^2+\xi^3+1+\frac{4}{27}t^{3n}=0.\]
The field automorphism defined by $t\rightarrow \sqrt[3n]{4/27}t$ brings this precisely to the curve $E^{1a}_{3n}$.
From this it follows that $E^{1a}_{3n}$ and $E^{5b}_n$ are $k$-equivalent.

We conclude that the curves $E^{1a}_{3n}$, $E^{1e}_{3n}$, $E^{5b}_{n}$, $E^{5c}_{n}$, $E^{5e}_{n}$, $E^{5g}_{n}$, $E^{5h}_{n}$ and $E^{5i}_{n}$ are all $k$-equivalent, and hence have the same rank.

\subsection{1b}
Here we show that $E^{1b}_{3n}$, $E^{4b}_{3n}$, $E^{4d}_{3n}$, $E^{7}_{3n}$ and $E^{8}_{2n}$ are $k$-equivalent. 
From this we can conclude that these families of curves have the same maximal rank.

The short Weierstrass form of $E^{7}_n$ is%can be reached by the the morphism given: $\eta=\sqrt{-4}(Y+\frac{1}{2})$ and $\xi=\sqrt[3]{-4}X$.
%This equation of the Weierstrass form is:
\[1+\frac{4}{\sqrt[3]{-4}}t^n\xi+\xi^3+\eta^2=0.\]
The field automorphism defined by $t\rightarrow (-4)^{2/(3n)}t$ brings this in the form described by $E^{1b}_n$.

Using ideas from \cite[Ch. 8 ]{cassels} we find that the curve $E^{4b}_n$ is isomorphic to the one given by:
\[\eta^2+\frac{t^n(i-1)}{2}\eta+\xi+\xi^3=0.\]
After a field automorphism this gives exactly the curve $E^7_n$.

The curve of the form $E^{4b}_n$ and $E^{4d}_n$ are isomorphic by the isomorphism $(X,Y)\rightarrow(\frac{1}{X},\frac{Y}{X^2})$.

The curve $E^8_{2n}$ is isomorphic to the curve given by:
\[\eta^2+1+\xi^3-\frac{t^{4n}}{4}\xi^4=0.\]
There is a field homomorphism bringing this precisely to $E^{4d}_{3n}$.

%This maps the curve of 8 into a subfamily of the curve of type 4d.
%We can conclude that $E^{1b}_{2n}$, $E^{4b}_{2n}$, $E^{4d}_{2n}$, $E^{7}_{2n}$ and $E^{8}_n$ are all $k$-equivalent.

\subsection{1c}
We will show that the curves $E^{1c}_{2n}$, $E^{1f}_n$ and $E^{6}_{2n}$ are $k$-equivalent.

The curve $E^{1f}_n$ is isomorphic to the curve given by: 
\[\eta^2+1+\xi^3-\frac{1}{4}t^{2n}\xi^2=0.\]
%As can be seen by the isomorphism $\eta=Y+\frac{1}{2}t^nX$ and $\xi=X$.
There is a field homomorphism bringing this curve precisely to $E^{1c}_{2n}$.
Likewise there is an isomorphism from $E^6_n$ to the curve given by %$\eta=2\sqrt{-1}Y+\sqrt{-1}$ and $\xi=\sqrt[3]{-4}X$ mapping the curve $E^6_n$ to the curve given by:
\[\eta^2+1+\xi^3+\sqrt[3]{-4}t^n\xi^2=0,\]
and there is a field homomorphism sending this curve to $E^{1c}_n$.

\subsection{1d}
We will show that the curves $E^{1d}_n$, $E^{3a}_n$, $E^{3b}_n$, $E^{3c}_n$, $E^{5a}_n$, $E^{5d}_n$ and  $E^{5j}_n$ are isomorphic.

Permuting $X$, $Y$, $Z$ gives isomorphisms between the curves $E^{5a}_n$, $E^{5d}_n$ and $E^{5j}_n$.

Likewise the curves $E^{3a}_n$, $E^{3b}_n$ and $E^{3c}_n$ are isomorphic by the morphisms:
$(X,Y)\rightarrow (X,(1+t^n)Y)$ from $E^{3c}_n$ to $E^{3b}_n$, and $(X,Y)\rightarrow ((1+t^n)X,(1+t^n)Y)$ from $E^{3b}_n$ to $E^{3a}_n$.

There is a isomorphism from $E^{3b}_n$ to $E^{1d}_n$ given by $(X,Y)\rightarrow(\sqrt[3]{-4}X,\sqrt{-1}(2Y+1))$.

Using ideas from \cite[Ch. 8 ]{cassels} we find that the curve $E^{5a}_n$ is isomorphic to $E^{1d}_n$.%hardest part left out

\subsection{2a}
We will show that the curves $E^{2a}_{2n}$, $E^{2c}_{2n}$, $E^{4a}_{2n}$, $E^{4e}_{2n}$, $E^{4f}_n$ and $E^{4h}_n$ are $k$-equivalent.

The curves $E^{2a}_n$ and $E^{2c}_n$ are isomorphic and the isomorphism from $E^{2a}_n$ to $E^{2c}_n$ is given by $(X,Y)\rightarrow (\frac{1}{X},\frac{Y}{X^2})$.
The curves $E^{4a}_n$ and $E^{4e}_n$ are also isomorphic, with isomorphism given by $(X,Y)\rightarrow (\frac{1}{X},\frac{Y}{X^2})$.
The curves $E^{4f}_n$ and $E^{4h}_n$ are likewise isomorphic, with isomorphism given by $(X,Y)\rightarrow (\frac{1}{X},\frac{Y}{X^2})$.

Using \cite[Ch. 8 ]{cassels} we bring $E^{4a}_n$ in short Weierstrass. This gives the curve $E^{2a}_n$.
So these two forms are also isomorphic.

The curve $E^{4f}_n$ can by taking $(X,Y)\rightarrow(X,Y+\frac{1}{2}t^n)$ be brought to the form:
\[\xi^4+\eta^2+1-\frac{1}{4}t^{2n}=0.\]
There is a field homomorphism bringing this precisely $E^{4a}_{2n}$.

\subsection{2b}

We show here that the curves $E^{2b}_n$, $E^{2d}_n$, $E^{4c}_{n}$, $E^{4g}_{n}$ and $E^{9}_n$ are $k$-isogenous and hence have the same rank.

The curve $E^{2d}_n$ is isomorphic to the curve defined by:
\[\eta^2+\xi^3+\xi^2+16t^n\xi=0.\]
%The morphism is given by $\eta=8\sqrt{-1}Y+4\sqrt{-1}X$ and $\xi=-4X$.
There is a field homomorphism sending this curve to $E^{2b}_n$, so $E^{2d}_n$ and $E^{2b}_n$ are $k$-equivalent.

There is an isogeny from $E^{4g}_{n}$ to $E^{2d}_n$ given by:
$(X,Y)\rightarrow (X^2,XY)$, so these two curves are isogenous.

%The curves of the forms 4c, 4g and 9 are all birational to each other.
There is an isomorphism from $E^{4g}_n$ to the curve given by
\[\eta^2+\xi^4+\xi^2+16t^n=0.\] 
%This isomorphism is given by $\xi=X$ and $\eta=Y+\frac{1}{2}t^nX$. 
There is a field homomorphism sending this curve to $E^{4c}_{n}$.

Likewise there is is an isomorphism from $E^{9}_n$ to the curve given by
\[\eta^2-\frac{1}{4}t^n+\xi^2+\xi^4=0.\]
%This isomorphism is given by $\xi=\frac{1}{\sqrt[4]{-4}}X$ and $\eta=Y+\frac{1}{2}X^2$.
There is a field homomorphism sending this curve to $E^{4c}_{n}$.

\subsection{2e}
The curves $E^{2e}_n$ and $E^{4i}_n$ are isomorphic.
An isomorphism from $E^{4i}_n$ to $E^{2e}_n$ is given by:
$(X,Y)\rightarrow (\frac{\sqrt[4]{-1}-X}{\sqrt[4]{-1}+X}, \frac{\sqrt{-2}Y}{(X+\sqrt[4]{-1})^2})$.

\subsection{3d}
There is an isogeny from $E^{5f}_n$ to $E^{3d}_n$, given by $(X,Y)\rightarrow(XY,Y^3)$.
Hence the curves $E^{5f}_n$ and $E^{3d}_n$ are isogenous.

\subsection{10}
There is also an isogeny from $E^{12}_n$ to $E^{10}_n$, given by $(X,Y)\rightarrow(XY^{-1},XY)$, hence the curves $E^{10}_n$ and $E^{12}_n$ are isogenous.

\section{Calculation of the ranks.}

In Appendix A we found a number of families that for various reasons have the same maximal rank. 
Here we will use Shioda's method to calculate these ranks.

\subsection{1a}

From Appendix A we know that the curves $E^{1e}_{3n}$, $E^{5b}_n$, $E^{5c}_n$, $E^{5e}_n$, $E^{5g}_n$, $E^{5h}_n$ and $E^{5i}_n$ are all $k$-isomorphic to a curve $E^{1a}_{3n}$.
In particular that means that these families all have the same maximal rank.
We will give the computation of the maximal rank of a curves in the family $E^{1a}_n$.

The maximal rank of this family was already computed by Shioda \cite{shioda2} in 1992. 
More details about this curve were computed in \cite{usui} and \cite{top} in 2000. 
We will repeat their results here out of a sense of completeness.
The curve $E^{1a}_n$ is defined by:
\[f=1+t^n+X^3+Y^2=0.\]
This is already in short Weierstrass form, with
\[\Delta=-432(1+t^n)^2, \mbox{ and }j=0.\]
The corresponding elliptic surface has a smooth fibre at $t=\infty$ precisely when $6|n$.
For the rest of this calculation we will assume $6|n$.
In this case the surface has precisely $n$ singular fibres of type \Rmnum{2}.
This gives $\rho_{triv}=2$.
We have that the second Betti number $h^2=2n-2$

Shioda's method can now be used to find $\lambda$.
We find the vectors generating $L$:
\[v_1=(-\frac{1}{3},0,\frac{1}{3},0),\;v_2=(-\frac{1}{2},0,0,\frac{1}{2}),\;v_3=(\frac{1}{n},-\frac{1}{n},0,0).\]

It can easily be shown that $iv_3,iv_3+v_1,iv_3+2v_1,iv_3+v_2\not\in \Lambda$.
%For $iv_3+v_1+v_2$ and $iv_3+2v_1+v_2$ it is not clear.
We will now determine when $iv_3+v_1+v_2$ is an element of $\Lambda$.
Write $t(iv_3+v_1+v_2)=(\frac{jt}{m},-\frac{jt}{m}-\frac{t5}{6},\frac{t}{3},\frac{t}{2})$.
Here $\frac{j}{m}=\frac{i}{n}-\frac{5}{6}$.
Compute: 
\[\fr{\frac{tj}{m}}+\fr{-\frac{tj}{m}-\frac{t5}{6}}+\fr{\frac{t}{3}}+\fr{\frac{t}{2}}=\left\{
\begin{array}{ll}
1& \mbox{ if } t\equiv 1\mod 6 \mbox{ and }  \fr{\frac{tj}{m}}<\frac{1}{6}\\
2& \mbox{ if } t\equiv 1\mod 6 \mbox{ and }  \fr{\frac{tj}{m}}>\frac{1}{6}\\
2& \mbox{ if } t\equiv 5\mod 6 \mbox{ and }  \fr{\frac{tj}{m}}<\frac{5}{6}\\
3& \mbox{ if } t\equiv 5\mod 6 \mbox{ and }  \fr{\frac{tj}{m}}>\frac{5}{6}\\
\end{array}
\right.
\]
This means that $iv_3+v_1+v_2\not\in \Lambda$ precisely when either $m\leq 6$ or $m>6$ and there does not exist a $j'\equiv 2\mod3$ such that $(j',m)=1$ and $\fr{\frac{j'}{m}}<\frac{1}{6}$.
The only $m$ for which this happens are $1, 2, 3, 4, 5, 6, 9, 12, 18, 24, 30, 60$.
For elements $iv_3+2v_1+v_2\in L$ we get a similar result.
This gives that $\lambda$ is at least $2n-72$, with equallity if $360|n$.

We can now compute the maximal rank:
\[r=h^{2}-\lambda-\rho_{triv}=2n-2-(2n-72)-2=68.\]

\subsection{1b}
In Appendix A we found that the curves $E^{1b}_{2n}$, $E^{4b}_{2n}$, $E^{4d}_{2n}$, $E^{7}_{2n}$ and $E^{8}_n$ are all $k$-equivalent.
This means that the maximal rank of these families of elliptic curves will be the same.
We will give the details of the computation of the maximal rank for the family $E^{1b}_n$.
Note that this example has already been treated by Shioda in \cite{shioda1}.

A maximal curve will be of the form:
\[f=1+t^nX+X^3+Y^2=0.\]

This is in short Weierstrass form so we can easily compute
\[\Delta=-64t^{3n}-432, \mbox{ and } j=1728\frac{4t^{3n}}{4t^{3n}+27}.\]
From this point we will assume $4|n$ so that the corresponding elliptic surface has a smooth fibre at $t=\infty$.
The surface has $3n$ fibres of type \Rmnum{1}$_1$ and no other singular fibres.
This gives $\rho_{triv}=2$.
The second Betti number can also be determined as $h^2=3n-2$.

We use Shioda's method to find $\lambda$.
Homogenizing $f$ gives:
\[Z^{n+1}+T^nX+X^3Z^{n-2}+Y^2Z^{n-1}=0.\]
From here we can compute the vectors generating $L$:
\[v_1=(-\frac{1}{2},0,0,\frac{1}{2}),\;v_2=(\frac{2}{3n},-\frac{1}{n},\frac{1}{3n},0).\]

We easily find that $iv_2\not\in \Lambda$.
For $iv_2+v_1$ we 
write $t(iv_2+v_1)=
(2\frac{jt}{m},-3\frac{jt}{m}-\frac{3t}{4},\frac{jt}{m}+\frac{t}{4},\frac{t}{2})$,
where $\frac{j}{m}=\frac{i}{3n}-\frac{1}{4}$.
Now compute 
\[\fr{2\frac{jt}{m}}+\fr{-3\frac{jt}{m}-\frac{3t}{4}}+\fr{\frac{jt}{m}+\frac{t}{4}}+\fr{\frac{t}{2}}=\left\{
\begin{array}{ll}
1& \mbox{ if } t\equiv 1\mod 4 \mbox{ and }  \fr{\frac{tj}{m}}<\frac{1}{12}\\
2& \mbox{ if } t\equiv 1\mod 4 \mbox{ and }  \frac{1}{12}<\fr{\frac{tj}{m}}<\frac{5}{12}\\
3& \mbox{ if } t\equiv 1\mod 4 \mbox{ and }  \frac{5}{12}<\fr{\frac{tj}{m}}<\frac{1}{2}\\
2& \mbox{ if } t\equiv 1\mod 4 \mbox{ and }  \fr{\frac{tj}{m}}>\frac{1}{2}\\
2& \mbox{ if } t\equiv 3\mod 4 \mbox{ and }  \fr{\frac{tj}{m}}<\frac{1}{2}\\
1& \mbox{ if } t\equiv 3\mod 4 \mbox{ and }  \frac{1}{2}<\fr{\frac{tj}{m}}<\frac{7}{12}\\
2& \mbox{ if } t\equiv 3\mod 4 \mbox{ and }  \frac{7}{12}<\fr{\frac{tj}{m}}<\frac{11}{12}\\
3& \mbox{ if } t\equiv 3\mod 4 \mbox{ and }  \fr{\frac{tj}{m}}>\frac{11}{12}\\
\end{array}
\right.
\]

From this we can see that that $iv_2+v_1\not\in \Lambda$ precisely when either $m\in \{1, 2, 3, 4, 5, 6, 8, 10, 12\}$ or $j\equiv 3\mod 4$ and $m \in \{20, 28, 36, 60, 84\} $.
This gives that $\lambda$ is at least $3n-60$, with equality if $840|n$.

Combining the results gives (for $840|n$) that:
\[r=h^{2}-\lambda-\rho_{triv}=3n-2-2-(3n-60)=56.\]

\subsection{1c}
Earlier we proved that the curves $E^{1c}_{2n}$, $E^{1f}_n$ and $E^{6}_{2n}$ are $k$-equivalent.
To find the maximal rank of these families of elliptic curves it suffices to find the maximal rank of the family of curves $E^{1c}_n$.

The curve $E^{1c}$ is given by:
\[f=1+t^nX^2+X^3+Y^2=0.\]
Bringing this to short Weierstrass gives:
\[\eta^2+\xi^3-\frac{t^{2n}}{3}\xi+1+\frac{2}{27}t^{3n}=0.\]
and we find
\[\Delta=-64t^{3n}-432, \mbox{ and }j=-\frac{256t^{6n}}{4t^{3n}+27}.\]
For the rest of this calculation we will assume that $2|n$, so that the only singular fibres of the corresponding elliptic surface are $3n$ fibres of type \Rmnum{1}$_1$ and one of type \Rmnum{1}$_{3n}$.
From here we see that $\rho_{triv}=1+3n$.
We can also compute the second Betti number $h^2=6n-2$

We use Shioda's method to find $\lambda$.
Homogenizing $f$ gives:
\[Z^{n+2}+T^nX^2+X^3Z^{n-1}+Y^2Z^{n}=0.\]
From here we determine the generators of $L$:
\[v_1=(-\frac{1}{2},0,0,\frac{1}{2}),\;v_2=(\frac{1}{3n},-\frac{1}{n},\frac{2}{3n},0).\]

We trivially find that $iv_2\not\in \Lambda$.
For $iv_2+v_1$ we write $t(iv_2+v_1)=
(\frac{jt}{m},-3\frac{jt}{m}+\frac{t}{2},2\frac{jt}{m},\frac{t}{2})$.
Here $\frac{j}{m}=\frac{i}{3n}-\frac{1}{2}$.
We compute
\[\fr{\frac{jt}{m}}+\fr{-3\frac{jt}{m}+\frac{t}{2}}+\fr{2\frac{jt}{m}}+\fr{\frac{t}{2}}=\left\{
\begin{array}{ll}
1& \mbox{ if } \fr{\frac{tj}{m}}<\frac{1}{6}\\
2& \mbox{ if } \frac{1}{6}<\fr{\frac{tj}{m}}<\frac{5}{6}\\
3& \mbox{ if } \frac{5}{6}<\fr{\frac{tj}{m}}\\
\end{array}
\right.
\]

This means that $iv_2+v_1\not\in \Lambda$ precisely when $m\in \{1, 2, 3, 4, 5, 6\}$.
This gives that $\lambda$ is at least $3n-12$, and equality holds if $20|n$.

Combining the results gives, for $20|n$ that the rank is:
\[r=h^{2}-\lambda-\rho_{triv}=6n-2-(3n+1)-(3n-12)=9.\]

\subsection{1d}
We already found that the curves $E^{1d}_n$, $E^{3a}_n$, $E^{3b}_n$, $E^{3c}_n$, $E^{5a}_n$, $E^{5d}_n$ and  $E^{5j}_n$ are $k$-invariant.
In the example we already computed the maximal rank of the curve $E^{1d}_n$ and found $r=18$.

\subsection{1g}

We will now compute the maximal rank of the family of curves $E^{1g}_n$.
The curve $E^{1g}_n$ is given by:
\[f=1+X^3+(1+t^n)Y^2=0.\]
Bringing this is in short Weierstrass form gives:
\[\eta^2+\xi^3+(1+t^n)^3=0.\]
And we can compute:
\[\Delta=-432(t^{n}+1)^6,\mbox{ and }j=0.\]
For the rest of the calculation we will assume $2|n$.
In that case the corresponding elliptic surface has exactly $n$ singular fibres of type \Rmnum{1}$_0^*$ and no other singular fibres.
From this we find $\rho_{triv}=2+4n$.
The second Betti number can also be found $h^2=6n-2$.

We use Shioda's method to find $\lambda$.
Homogenizing $f$ gives:
\[Z^{n+2}+X^3Z^{n-1}+Y^2Z^{n}+Y^2T^n=0.\]
This gives generators for $L$:
\[v_1=(-\frac{1}{3},\frac{1}{3},0,0),\;v_2=(-\frac{1}{2},0,\frac{1}{2},0),\;v_3=(0,0,\frac{1}{n},-\frac{1}{n}).\]

It can easilly be seen that $iv_3,v_1+iv_3,2v_1+iv_3, v_2+iv_3\not\in \Lambda$.
%For $v_1+v_2+iv_3$ and $2v_1+v_2+iv_3$ it is not clear.
We will now determine when $v_1+v_2+iv_3 \in \Lambda$.
Write $t(v_1+v_2+iv_3)=
(-\frac{5t}{6},\frac{t}{3},\frac{t}{2}+\frac{jt}{m},-\frac{jt}{m})$,
where $j$ and $m$ are minimal such that $j/m=i/n$.
Compute 
\[\fr{-\frac{5t}{6}}+\fr{\frac{t}{3}}+\fr{\frac{t}{2}+\frac{jt}{m}}+\fr{-\frac{jt}{m}}=\left\{
\begin{array}{ll}
2& \mbox{ if } t\equiv 1 \mod 6 \mbox{ and }\fr{\frac{tj}{m}}<\frac{1}{2}\\
1& \mbox{ if } t\equiv 1 \mod 6 \mbox{ and }\fr{\frac{tj}{m}}>\frac{1}{2}\\
3& \mbox{ if } t\equiv 5 \mod 6 \mbox{ and }\fr{\frac{tj}{m}}<\frac{1}{2}\\
2& \mbox{ if } t\equiv 5 \mod 6 \mbox{ and }\fr{\frac{tj}{m}}>\frac{1}{2}\\
\end{array}
\right.
\]

This means that $v_1+v_2+iv_3\not\in \Lambda$ precisely when $m\in \{1, 2\} $ or when $m\in \{3, 6\}$ and $j\equiv 1\mod 3$. For $2v_1+v_2+iv_3$ we get a similar result.
This gives that $\lambda$ is at least $2n-8$, and equality holds if $6|n$.

Combining this gives if $6|n$ the rank:
\[r=h^{2}-\lambda-\rho_{triv}=6n-2-(2n-8)-(4n+2)=4.\]

\subsection{2a}
Previously we found that the curves $E^{2a}_{2n}$, $E^{2c}_{2n}$, $E^{4a}_{2n}$, $E^{4e}_{2n}$, $E^{4f}_n$ and $E^{4h}_n$ are all $k$-equivalent.
We will only have to compute the maximal rank of the family of curves $E^{2a}_n$.
The curve $E^{2a}_{n}$ is defined by:
\[f=(1+t^n)X+X^3+Y^2=0.\]
This is already in short Weierstrass form so we can easily compute.
\[\Delta=-64(t^{n}+1)^3, \mbox{ and } j=1728.\]
From here we will assume $4|n$.
In this case the corresponding elliptic surfaces has $n$ singular fibres of type \Rmnum{3} and no other singular fibres.
This gives $\rho_{triv}=2+n$.
The second Betti number will be $h^2=3n-2$.

We use Shioda's method to find $\lambda$.
Homogenizing $f$ gives:
\[XZ^{n}+XT^n+X^3Z^{n-2}+Y^2Z^{n-1}=0.\]
From this we can compute generators for $L$:
\[v_1=(-\frac{3}{4},0,\frac{1}{4},\frac{1}{2}),\;v_2=(\frac{1}{n},-\frac{1}{n},0,0).\]

It can be easily seen that $iv_2,2v_1+iv_2\not\in \Lambda$.
%For $v_1+iv_2$ and $3v_1+iv_2$ it is not clear.
We now have to determine when $v_1+iv_2\in \Lambda$.
Write $t(v_1+iv_2)=
(\frac{jt}{m}-\frac{3t}{4},-\frac{jt}{m},\frac{t}{4},\frac{t}{2})$,
where $j,m$ are minimal such that $j/m=i/n$.
Now compute 
\[\fr{\frac{jt}{m}-\frac{3t}{4}}+\fr{-\frac{jt}{m}}+\fr{\frac{t}{4}}+\fr{\frac{t}{2}}=\left\{
\begin{array}{ll}
2& \mbox{ if } t\equiv 1 \mod 4 \mbox{ and }\fr{\frac{tj}{m}}<\frac{3}{4}\\
1& \mbox{ if } t\equiv 1 \mod 4 \mbox{ and }\fr{\frac{tj}{m}}>\frac{3}{4}\\
3& \mbox{ if } t\equiv 3 \mod 4 \mbox{ and }\fr{\frac{tj}{m}}<\frac{1}{4}\\
2& \mbox{ if } t\equiv 3 \mod 4 \mbox{ and }\fr{\frac{tj}{m}}>\frac{1}{4}\\
\end{array}
\right.
\]
This means that $v_1+iv_2\not\in \Lambda$ precisely when $m\in \{1,2,3,4\} $ or when $m\in \{8, 12, 24\}$ and $j\equiv 1\mod 4$. For $3v_1+iv_2$ we get a similar result.
This gives that $\lambda$ is at least $2n-28$, and equality holds if $24|n$.

Combining the we find that if $24|n$ the rank:
\[r=h^{2}-\lambda-\rho_{triv}=3n-2-(2n-28)-(n+2)=24.\]

\subsection{2b}
Previously we proved that the curves $E^{2b}_n$, $E^{2d}_n$, $E^{4c}_{n}$, $E^{4g}_{n}$ and $E^{9}_n$ are $k$-isogenous, and as such have the same rank.
We will compute the maximal rank of the family of curves $E^{2b}_n$.
The curve $E^{2b}_n$ is given by
\[f=t^nX+X^2+X^3+Y^2=0.\]
In short Weierstrass form this curve is given by:
\[\eta^2+\xi^3+(t^n-\frac{1}{3})\xi+\frac{2}{27}-\frac{t^n}{3}=0.\]
From this we can compute:
\[\Delta=-64t^{3n}+16t^{2n},\mbox{ and }j=256\frac{(3t^n-1)^3}{4t^{3n}-t^{2n}}.\]
From here on we will assume that $4|n$.
In this case the corresponding elliptic surface has $n$ singular fibres of type \Rmnum{1}$_1$, one singular fibre of type \Rmnum{1}$_{2n}$ and no other singular fibres.
This gives $\rho_{triv}=2n+1$.
The second Betti number is given by $h^2=3n-2$.

We use Shioda's method to find $\lambda$.
Homogenizing $f$ gives:
\[XT^n+X^2Z^{n-1}+X^3Z^{n-2}+Y^2Z^{n-1}=0.\]
This gives generators for $L$:
\[v_1=(0,-\frac{1}{2},0,\frac{1}{2}),\;v_2=(-\frac{1}{n},\frac{2}{n},-\frac{1}{n},0).\]

It is easily seen that $iv_2\not\in \Lambda$.
For $v_1+iv_2$ we write $t(v_1+iv_2)=
(-\frac{jt}{m},2\frac{jt}{m}-\frac{t}{2},-\frac{jt}{m},\frac{t}{2})$.
Here $j,m$ are minimal such that $j/m=i/n$.
We can compute:
\[\fr{-\frac{jt}{m}}+\fr{2\frac{jt}{m}-\frac{t}{2}}+\fr{-\frac{jt}{m}}+\fr{\frac{t}{2}}=\left\{
\begin{array}{ll}
3& \mbox{ if } \fr{\frac{tj}{m}}<\frac{1}{4}\\
2& \mbox{ if } \frac{1}{4}<\fr{\frac{tj}{m}}<\frac{3}{4}\\
1& \mbox{ if } \frac{3}{4}<\fr{\frac{tj}{m}}\\
\end{array}
\right.
\]
This means that $v_1+iv_2\not\in \Lambda$ precisely when $m\in \{1,2,3,4\} $.
This gives that $\lambda$ is at least $n-6$, and equality holds if $12|n$.

It follows that if $12|n$ the rank is:
\[r=h^{2}-\lambda-\rho_{triv}=3n-2-(n-6)-(2n+1)=3.\]

\subsection{2e}
The curves $E^{2e}_n$ and $E^{4i}_n$ are isomorphic, as such they have the same rank.
We will compute the maximal rank of the family $E^{2e}_n$.
The curve $E^{2e}_n$ is defined by:
\[f=X+X^3+(1+t^n)Y^2=0.\]
In short Weierstrass form this gives:
%By taking $\xi=(1+t^n)X$ and $\eta=(1+t^n)^2Y$.
\[\eta^2+\xi^3+(t^n+1)^2\xi=0.\]
For this curve we can compute:
\[\Delta=-64(1+t^n)^6, \mbox{ and }j=1728.\]
We will from here on assume that $2|n$.
In that case the corresponding elliptic surface has $n$ singular fibres of type \Rmnum{1}$_0^*$ and no other singular fibres.
This gives $\rho_{triv}=4n+2$.
The second Betti number can be determined $h^2=6n-2$.

We use Shioda's method to find $\lambda$.
Homogenizing $f$ gives:
\[XZ^{n+1}+X^3Z^{n-1}+Y^2Z^{n}+Y^2T^n=0.\]
This gives as generators for $L$:
\[v_1=(\frac{1}{4},\frac{1}{4},\frac{1}{2},0),\;v_2=(0,0,\frac{1}{n},-\frac{1}{n}).\]

It is easy to see that $iv_2, 2v_1+iv_2\not\in \Lambda$.
Write $t(v_1+iv_2)=
(\frac{t}{4},\frac{t}{4},\frac{jt}{m}+\frac{t}{2},-\frac{jt}{m})$, where $j,m$ are minimal such that $j/m=i/n$.
We can compute 
\[\fr{\frac{t}{4}}+\fr{\frac{t}{4}}+\fr{\frac{jt}{m}+\frac{t}{2}}+\fr{-\frac{jt}{m}}=\left\{
\begin{array}{ll}
2& \mbox{ if } t\equiv 1 \mod 4 \mbox{ and } \fr{\frac{tj}{m}}<\frac{1}{2}\\
1& \mbox{ if } t\equiv 1 \mod 4 \mbox{ and } \fr{\frac{tj}{m}}>\frac{1}{2}\\
3& \mbox{ if } t\equiv 3 \mod 4 \mbox{ and } \frac{1}{2}>\fr{\frac{tj}{m}}\\
2& \mbox{ if } t\equiv 3 \mod 4 \mbox{ and } \frac{1}{2}<\fr{\frac{tj}{m}}\\
\end{array}
\right.
\]
This means that $v_1+iv_2\not\in \Lambda$ precisely when $m\in \{1,2\} $ or when $m\in\{4, 12\}$ and $j\equiv 1\mod 4$. 

For $3v_1+iv_2$ we get a similar result. 
Now it follows that $\lambda$ is at least $2n-10$, and equality holds if $12|n$.
We conclude that if $12|n$ the rank of our curve is:
\[r=h^{2}-\lambda-\rho_{triv}=6n-2-(2n-10)-(4n+2)=6.\]

\subsection{3d}
The curves $E^{5f}_n$ and $E^{3d}_n$ are isogenous, as such they have the same rank.
We will now compute the maximal rank of the family of curves $E^{3d}_n$.
The curve $E^{3d}_n$ is given by:
\[f=Y+t^nXY+X^3+Y^2=0.\]
Bringing this in short Weierstrass form gives:
\[\eta^2+\xi^3-(\frac{1}{48}t^{4n}+\frac{1}{2}t^n)\xi-\frac{1}{4}-\frac{1}{24}t^{3n}-\frac{1}{864}t^{6n}=0.\]
For this form we compute the invariants:
\[\Delta=-(t^{3n}+27), \mbox{ and } j=-\frac{(t^{4n}+24t^n)^3}{t^{3n}+27}.\]
This has $3n$ singular fibres of type \Rmnum{1}$_1$, one singular fibre of type \Rmnum{1}$_{9n} $.
From this we can compute $\rho_{triv}=9n+1$.
The second Betti number can also be determined $h^2=12n-2$.%, when $2|n$.

We use Shioda's method to find $\lambda$.
Homogenizing $f$ gives:
\[YZ^{n+1}+T^nXY+X^3Z^{n-1}+Y^2Z^{n}=0.\]
Using this we can find that $L$ is generated by:
\[v_1=(\frac{1}{3n},-\frac{1}{n},\frac{1}{3n},\frac{1}{3n}).\]
It has to determine whether $iv_1\in \Lambda$ or not.
Write $t(iv_1)=
(\frac{jt}{m},-3\frac{jt}{m},\frac{jt}{m},\frac{jt}{m})$.
Here $j,m$ are minimal such that $j/m=i/3n$.
We can compute
\[\fr{\frac{jt}{m}}+\fr{-3\frac{jt}{m}}+\fr{\frac{jt}{m}}+\fr{\frac{jt}{m}}=\left\{
\begin{array}{ll}
1& \mbox{ if } \fr{\frac{tj}{m}}<\frac{1}{3}\\
2& \mbox{ if } \frac{1}{3}<\fr{\frac{tj}{m}}<\frac{1}{3}\\
3& \mbox{ if } \frac{2}{3}<\fr{\frac{tj}{m}}\\
\end{array}
\right.
\]
This means that $iv_1\not\in \Lambda$ precisely when $m\in \{1,2,3\}$. 
This gives that $\lambda$ is at least $3n-4$, and equality holds if $n$ is even.
We conclude that for even $n$ the rank will be:
\[r=h^{2}-\lambda-\rho_{triv}=12n-2-(3n-4)-(9n+1)=1.\]

\subsection{11}
We will determine the maximal rank of curves of the family $E^{11}_n$.
The curve $E^{11}_n$ is given by.
\[f=t^n+XY^2+X^3+Y^2=0.\]
We will for our calculation assume that $6|n$.
In this case we use ideas from \cite[Ch. 8 ]{cassels} to find that the short Weierstrass form becomes:
%Done by first making isomorphic to Y^2+t^n+t^nX+X^3+X^4 and the Cassel's trick.
%By taking $\xi=X/Y+\frac{1}{3}(1+t^n)$ and $\eta=Y+1$.
%After that dividing by t^{a}
\[\eta^2+\xi^3-3t^{n/3}\xi+1+t^{n}=0.\]
This has the following invariants:
\[\Delta=-432(1-t^{n})^2\mbox{ and }j=-\frac{6912t^{n}}{(1-t^{n})^2}.\] 
The corresponding surface has exactly $n$ singular fibres all of type \Rmnum{1}$_2$.
It follows that $\rho_{triv}=n+2$.
The second Betti number in this case is $h^2=2n-2$.

We use Shioda's method to find $\lambda$.
Homogenizing $f$ gives:
\[T^n+XY^2Z^{n-3}+X^3Z^{n-3}+Y^2Z^{n-2}=0.\]
From here we compute generators for $L$:
\[v_1=(0,\frac{1}{2},\frac{1}{2},0),\;v_2=(-\frac{1}{n},-\frac{3}{n},\frac{1}{n},\frac{3}{n}).\]

It is easily seen that $iv_2, \not\in \Lambda$.
For $v_1+iv_2$ we write $t(v_1+iv_2)=
(-\frac{jt}{m}-\frac{t}{6},-3\frac{jt}{m},\frac{2t}{3}+\frac{jt}{m},\frac{t}{2}+3\frac{jt}{m})$.
Here $j,m$ are minimal such that $j/m=i/n-1/6$.
We can compute:
\[\fr{-\frac{jt}{m}-\frac{t}{6}}+\fr{-3\frac{jt}{m}}+\fr{\frac{2t}{3}+\frac{jt}{m}}+\fr{\frac{t}{2}+3\frac{jt}{m}}=\left\{
\begin{array}{ll}
3& \mbox{ if } t\equiv 1\mod 6 \mbox{ and } \fr{\frac{tj}{m}}<\frac{1}{6}\\
2& \mbox{ if } t\equiv 1\mod 6 \mbox{ and } \frac{1}{6}<\fr{\frac{tj}{m}}<\frac{1}{2}\\
1& \mbox{ if } t\equiv 1\mod 6 \mbox{ and } \frac{1}{2}<\fr{\frac{tj}{m}}<\frac{2}{3}\\
2& \mbox{ if } t\equiv 1\mod 6 \mbox{ and } \frac{2}{3}<\fr{\frac{tj}{m}}\\
2& \mbox{ if } t\equiv 5\mod 6 \mbox{ and } \fr{\frac{tj}{m}}<\frac{1}{3}\\
3& \mbox{ if } t\equiv 5\mod 6 \mbox{ and } \frac{1}{3}<\fr{\frac{tj}{m}}<\frac{1}{2}\\
2& \mbox{ if } t\equiv 5\mod 6 \mbox{ and } \frac{1}{2}<\fr{\frac{tj}{m}}<\frac{5}{6}\\
1& \mbox{ if } t\equiv 5\mod 6 \mbox{ and } \frac{5}{6}<\fr{\frac{tj}{m}}\\
\end{array}
\right.
\]
This means that $v_1+iv_2\not\in \Lambda$ precisely when $m\in \{1,2,3 ,4,6 \}$. or when $m\in\{12, 24, 60\}$ and $j\equiv 2\mod 3$. %For $3v_1+iv_2$ we get a similar result. 
This gives that $\lambda$ is at least $n-22$, and equality holds if $120|n$.
Combining these results gives if $120|n$:
\[r=h^{2}-\lambda-\rho_{triv}=2n-2-(n-22)-(n+2)=18.\]

\subsection{12}
Finally we will look at the family $E^{12}_n$.
This family consists of elliptic curves in the Edwards form.
The curve $E^{12}_n$ is given by:
\[f=t^n+X^2+Y^2+X^2Y^2=0.\]
%Going to short Weierstrass form is not so easy.
%First take %%%%%%%%%%
Using ideas from \cite[Ch. 8 ]{cassels} we reduce this curve to short Weierstrass form. 
This  gives:
%Done by first making isomorphic to Y^2+t^n+t^nX+X^3+X^4 and the Cassel's trick.
%By taking $\xi=X/Y+\frac{1}{3}(1+t^n)$ and $\eta=Y+1$.
%After that dividing by t^{a}
\[\eta^2+\xi^3-(\frac{1}{3}+\frac{14}{3}t^{n}+\frac{1}{3}t^{2n})\xi+\frac{2}{27}-\frac{22}{9}t^{n}-\frac{22}{9}t^{2n}+\frac{2}{27}t^{3n}=0.\]
For this we compute the invariants:
\[\Delta=-256 t^n(1-t^n)^4, \mbox{ and }j=-16\frac{(1+14t^{n}+t^{2n})^3 }{t^n(1-t^n)^4}.\]
We will for the rest of this calculation assume that $n$ is even.
The surface now has $n$ singular fibres of type \Rmnum{1}$_4$, $2$ singular fibres of type \Rmnum{1}$_{n}$ and no other singular fibres.
This gives $\rho_{triv}=5n$.
From here on we will assume that $n$ is even.
Under this assumption the second Betti number  is $h^2=6n-2$.

We use Shioda's method to find $\lambda$.
Homogenizing $f$ gives:
\[T^n+X^2Z^{n-2}+Y^2Z^{n-2}+X^2Y^2Z^{n-4}=0.\]
This gives the following generators for $L$:
\[v_1=(0,0,\frac{1}{2},\frac{1}{2}),\;v_2=(0,\frac{1}{2},0,\frac{1}{2}),\;v_3=(-\frac{1}{n},\frac{1}{n},\frac{1}{n},-\frac{1}{n}).\]

It turns out that $iv_3, v_1+iv_3, v_2+iv_3 \not\in \Lambda$.
Write $t(v_1+v_2+iv_3)=
(-\frac{jt}{m},\frac{jt}{m}+\frac{t}{2},\frac{jt}{m}+\frac{t}{2},-\frac{jt}{m})$, where $j,m$ are minimal such that $j/m=i/n$.
Now compute 
\[\fr{-\frac{jt}{m}}+\fr{\frac{jt}{m}+\frac{t}{2}}+\fr{\frac{jt}{m}+\frac{t}{2}}+\fr{-\frac{jt}{m}}=\left\{
\begin{array}{ll}
1& \mbox{ if } \frac{1}{2}<\fr{\frac{tj}{m}}\\
3& \mbox{ if } \fr{\frac{tj}{m}}<\frac{1}{2}\\
\end{array}
\right.
\]
This means that $v_1+v_2+iv_3\not\in \Lambda$ precisely when $m\in \{1,2\}$.% or when $m\in{12, 24, 60}$ and $j\equiv 2\mod 4$. %For $3v_1+iv_2$ we get a similar result. 
This gives that $\lambda$ is at least $n-2$, and equality holds if $2|n$.

Combining the results gives for even $n$:
\[r=h^{2}-\lambda-\rho_{triv}=6n-2-(n-2)-5n=0.\]

\end{document}